\documentclass[11pt,reqno]{amsart}
\usepackage{a4wide}
\usepackage{euscript,amsmath,amssymb,amsbsy,mathabx}
\usepackage[arrow,matrix,curve]{xy}
\usepackage{array,longtable,graphicx,enumitem,url}

\usepackage[colorlinks=true,linkcolor=black,urlcolor=black,citecolor=black]{hyperref}

\newcounter{msct}[section]\renewcommand{\themsct}{\thesection.\arabic{msct}}
\newenvironment{m-theorem}{\vskip5pt\refstepcounter{msct}\trivlist \itemindent 0pt%
\item[\hskip\labelsep\bf Theorem~\themsct]\it\ignorespaces}{\endtrivlist\vskip3pt}
\newenvironment{m-proposition}{\vskip5pt\refstepcounter{msct}\trivlist \itemindent0pt%
\item[\hskip\labelsep\bf Proposition~\themsct]\it\ignorespaces}{\endtrivlist\vskip3pt}
\newenvironment{m-corollary}{\vskip5pt\refstepcounter{msct}\trivlist \itemindent 0pt%
\item[\hskip\labelsep\bf Corollary~\themsct]\it\ignorespaces}{\endtrivlist\vskip3pt}
\newenvironment{m-lemma}{\vskip5pt\refstepcounter{msct}\trivlist \itemindent 0pt%
\item[\hskip\labelsep\bf Lemma~\themsct]\it\ignorespaces}{\endtrivlist\vskip3pt}
\newenvironment{m-definition}{\vskip5pt\refstepcounter{msct}\trivlist \itemindent0pt%
\item[\hskip\labelsep\bf Definition~\themsct]\ignorespaces}{\endtrivlist\vskip5pt}
\newenvironment{m-notation}{\vskip5pt\refstepcounter{msct}\trivlist \itemindent0pt%
\item[\hskip\labelsep\bf Notation~\themsct]\ignorespaces}{\endtrivlist\vskip5pt}
\newenvironment{m-example}{\vskip5pt\refstepcounter{msct}\trivlist \itemindent0pt%
\item[\hskip\labelsep\bf Example~\themsct]\ignorespaces}{\endtrivlist\vskip5pt}
\newenvironment{m-remark}{\vskip5pt\refstepcounter{msct}\trivlist \itemindent0pt%
\item[\hskip\labelsep\bf Remark~\themsct]\ignorespaces}{\endtrivlist\vskip5pt}
\newenvironment{m-question}{\vskip3pt\refstepcounter{msct}\trivlist \itemindent0pt%
\item[\hskip\labelsep\bf Question.]\ignorespaces}{\endtrivlist\vskip5pt}
\newenvironment{thm-nono}[1]{\vskip5pt\trivlist \itemindent 0pt %
\item[\hskip\labelsep\bf Theorem~{\rm\mbox{#1}}]\it\ignorespaces}{\endtrivlist\vskip5pt}
\newenvironment{lm-nono}[1]{\vskip5pt\trivlist \itemindent0pt%
\item[\hskip\labelsep\bf Lemma~{\rm\mbox{#1}}]\it\ignorespaces}{\endtrivlist\vskip5pt}
\newenvironment{conj-nono}[1]{\vskip5pt\trivlist \itemindent0pt%
\item[\hskip\labelsep\bf Conjecture~{\rm\mbox{#1}}]\it\ignorespaces}{\endtrivlist\vskip5pt}
\newenvironment{m-thank}{\vskip5pt\trivlist \itemindent0pt%
\item[\hskip\labelsep\it Acknowledgments]\ignorespaces}{\endtrivlist\vskip5pt}
\newenvironment{m-proof}{\vskip2pt\trivlist \itemindent0pt%
\item[\hskip\labelsep\it Proof.]\ignorespaces}{\hfill$\Box$\endtrivlist\vskip5pt}%
\newenvironment{m-asmp}{\vskip5pt\trivlist \itemindent0pt%
\item[\hskip\labelsep\bf Assumption.]\ignorespaces}{\hfill\endtrivlist\vskip5pt}%

\newcounter{meqn}[section]\renewcommand{\themeqn}{\thesection.\arabic{meqn}}
\newenvironment{m-eqn}[1]{\vskip5pt\refstepcounter{meqn}%
\trivlist\itemindent0pt\item[]\ignorespaces%
\hfill $\displaystyle #1$\hfill\hbox{\rm(\themeqn)}}{\endtrivlist\vskip5pt}

\newcommand{\bibauth}[2]{\textrm{{#1}~{#2}}}
\newcommand{\bibtitl}[1]{\textit{#1}.}
\newcommand{\bibjnyp}[4]{\textrm{#1} \textbf{#2} (#3), #4.}
\newcommand{\bibinbook}[4]{In: \textrm{#1}\textrm{, #2}\textrm{, #3}\textrm{, #4}.}

\let\lar\longrightarrow
\let\hra\hookrightarrow
\let\mt\mapsto

\let\dashto\dashrightarrow

\let\euf\EuScript %use of package ``euscript'' required
\let\cal\mathcal
\let\mbb\mathbb

\DeclareFontFamily{OT1}{rsfs}{}
\DeclareFontShape{OT1}{rsfs}{n}{it}{<->rsfs10}{}
\DeclareMathAlphabet{\crl}{OT1}{rsfs}{n}{it}

\numberwithin{equation}{section}\numberwithin{figure}{section}

\newcommand\bone{{1\kern-0.57ex\rm l}}

\newcommand\beE{{{\euf E}}}

\newcommand\eF{{\euf F}}
\newcommand\beF{{{\euf F}}}

\newcommand\beG{{{\euf G}}}

\newcommand{\Gl}{{\mathop{\rm GL}\nolimits}}
\newcommand\eI{{\euf I}}
\newcommand{\eK}{{\euf K}}
\newcommand\eO{{\euf O}}
\newcommand{\eQ}{{\euf Q}}
\newcommand{\eS}{{\euf S}}
\newcommand{\eT}{{\euf T}}
\newcommand{\cU}{{\cal U}}
\newcommand\dbar{{\bar{\partial }}}
\newcommand{\ev}{\mathop{{\rm ev}}\nolimits}
\let\veps\varepsilon 
\newcommand\End{\mathop{\rm{End}}\nolimits}
\newcommand\cEnd{\operatorname{{\cal E}\kern-1pt\textit{nd}\kern1pt}}
\newcommand\Ext{\mathop{\rm{Ext}}\nolimits}
\newcommand\cExt{\operatorname{{\cal E}\kern-1pt\textit{xt}\kern1pt}}
\newcommand{\Hom}{\mathop{{\rm Hom}}\nolimits}
\newcommand{\cHom}{\mathop{{\cal H}om}\nolimits}
\newcommand{\kz}{{{\kappa}\!z}}
\let\lda\lambda
\newcommand{\Ker}{{\mathop{\rm Ker}\nolimits}}
\let\ges\geqslant
\let\les\leqslant 
\let\nit\noindent 
\newcommand{\pl}{{\mathop{{\rm p}_{l}}\nolimits}}
\newcommand{\pr}{{\mathop{{\rm p}_{r}}\nolimits}}
\newcommand{\ppp}{{\mbb P^3}}
\newcommand{\PGL}{{\mathop{\rm PGL}\nolimits}}
\newcommand{\EgI}{{\mathop{\mbox{$\cal E$\rm I}}\nolimits}}

\newcommand\ouset[3]{{\overset{#2}{\underset{#1}#3}}}
\newcommand{\rst}{{\!\upharpoonright}}
\let\si\sigma 
\let\surj\twoheadrightarrow
\let\srel\stackrel 
\let\tld\tilde 
\let\wtld\widetilde 
\newcommand\arr[1]{\mathrel{\vrule height3.25pt depth-2.75pt width #1\kern-1ex\hbox{$\rightarrow$}}}

\title[Rationality of the moduli space of instantons on $\mbb P^3$]%
{On the rationality of the moduli space of\\ instanton bundles on the projective 3-space}
\author{M.~Halic, R.~Tajarod}
\email{mihai.halic@gmail.com, roshan.tajarod@gmail.com}
\address{CRM, UMI 3457, Montr\'eal~H3T~1J4, Canada}
\subjclass[2010]{14M20, 14J60, 14D20, 14D21}
\keywords{rationality; moduli; instantons}
\begin{document}
\begin{abstract}
We prove the rationality and irreducibility of the moduli space of---what we call---the 
endomorphism-general instanton vector bundles of arbitrary rank on $\mbb P^3$. In 
particular, we deduce the rationality of the moduli spaces of rank-$2$ mathematical 
instantons. This problem was first studied by Hartshorne, Hirschowitz-Narasimhan in the 
late 1970s, and it has been reiterated within the framework of the ICM~2018.
\end{abstract}
\maketitle 

\section*{Introduction}
The interest in rank-$2$ instanton bundles on the three-dimensional projective space (with 
Chern classes $c_1=0,c_2=n$) has its origins in the articles of Atiyah \textit{et 
al.}~\cite{adhm,ahs,atiyah-ward}, Barth-Hulek~\cite{brth+hulk,brth-inst}, and 
Hartshorne~\cite{hart-vb-P3}, which in turn were motivated by work of 
't~Hooft~\cite{hooft} and Polyakov~\cite{plkv}. The geometry of their moduli spaces, such 
as the irreducibility and rationality, has been intensively investigated, especially 
during the past decade. So far, it is known that these quasi-projective varieties are rational for  
$n=2,4,5$, due to works of Hirschowitz-Narasimhan~\cite{hirs+nars}, 
Ellingsrud-Str{\o}mme~\cite{ellin+strm}, and Katsylo~\cite{kats}; for $n\ges 6$, the issue 
remained open, in spite of efforts of Tikhomirov~\textit{et al.} 
(cf.~\cite{tikh-odd,tikh-even}). Let us remark that the techniques used in these 
references \emph{are specific to the rank-two case}. Recently, within the framework of the 
ICM~2018~\cite{tikh-2018}, it was reaffirmed the importance and necessity of addressing 
the rationality of the moduli space of instanton bundles on $\ppp$.

The first author investigated~\cite{hlc-sheaves} the geometry of the moduli space of semi-stable vector bundles on $\mbb P^2$-bundles over $\mbb P^1$. Our goal in here is to deepen those results in the case of the three-dimensional projective space and to address the rationality issue mentioned above. 
Note that, in spite of the extensive literature on the rank $r=2$ case, there are few articles dealing with higher rank bundles on $\mbb P^3$. However, non-abelian gauge theories frequently appear in the physics literature, especially for the special unitary group $SU(r)$, $r\ges2$. Also, the Penrose transform which relates Hermitian vector bundles endowed with self-dual connections over the sphere $S^4$, on the one hand, to certain holomorphic vector bundles on $\ppp$, on the other hand, does apply in this more general setting~\cite{ahs}. For this reason we believe that our unified treatment of the \emph{arbitrary rank case}, not only two, provides additional interest to this article.  

Now let us justify the reason for restricting our attention to a particular class of stable bundles on $\ppp$. Indeed, already for $r=2$, the description of the components of the moduli spaces of all the stable vector bundles is involved; the observation goes back to Barth-Hulek~\cite[\S8]{brth+hulk}. This fact forces us to single out a class of better behaved vector bundles, those that we call $\cEnd$-general instantons; in short, this means that both the vector bundle and its endomorphisms satisfy the so-called instanton condition (cf. Definition~\ref{def:FonP3} below). A good feature is that the corresponding moduli space has the expected dimension, given by the Riemann-Roch formula, and it contains (in the case $r=2$) the special 't Hooft bundles studied by Hirschowitz and Narasimhan~\cite{hirs+nars,beo+fra}.

\begin{thm-nono}{}
The moduli space of $\cEnd$-general instanton vector bundles on $\mbb P^3$ (of rank $r$, with Chern classes $0,n,0$) is non-empty, irreducible of dimension $4rn-r^2+1$, and it is \underbar{\emph{rational}}. In particular, the statement holds for the moduli space of rank-$2$ mathematical instanton bundles on $\mbb P^3$.
\end{thm-nono}
Furthermore, we give a completely explicit description of the generic $\cEnd$-general instanton bundle; it is determined by its restriction either to a pair of planes or to a smooth quadric in $\ppp$ (cf. Theorem~\ref{thm:B}). 

We conclude this introduction with a brief survey of the article. Our techniques are based on cohomological computations but, compared to previous work, our approach to the problem is original and it contains novelties: 
\begin{itemize}[leftmargin=3ex]
\item[--] 
The analysis of the properties of instantons requires an in-depth understanding of  Barth and Hulek's construction~\cite{brth+hulk}. We show that their monad is obtained in a functorial way, from a `universal' diagram which is based on Beilinson's resolution of the diagonal in $\ppp\times\ppp$; this allows to control various homomorphisms in cohomology, which appear in the display of the monad. To our knowledge, this fact \emph{has never been observed before}.\smallskip 
\item[--] 
For studying the geometry of their moduli space, we restrict instantons on $\ppp$ to (a pair of) planes. The $\cEnd$-general condition ensures that this restriction map is \'etale, so we can apply the results of~\cite{hlc-sheaves} concerning the rationality of the moduli space of framed instantons on $\mbb P^2$.
\end{itemize}
We also highlight the brevity of our work, making the proofs easy to follow, which reflects the effectiveness of our methods: the arguments are not restricted only to the rank-$2$ case; we deduce the irreducibility of the relevant moduli space in a couple of lines (cf. Theorem~\ref{thm:irred}). 

%%%%%%%%%%%%%%%%%%%%%%%%%%%%%%%%
%%%%%%%%%%%%%%%%%%%%%%%%%%%%%%%%

\section{The framework}\label{sct:frame}

Throughout the article, we work over an algebraically closed field $\Bbbk$ of characteristic zero. 

\begin{m-definition}\label{def:FonP3} 
An \emph{instanton-like} vector bundle $\beF$ on $\mbb P^3$ of charge $n$ and rank $r$, with $n\ges r$, is defined by the following properties:
\begin{enumerate}[leftmargin=7ex]
\item[(i)] 
 $c_1(\beF) = 0, c_2(\beF)=n, c_3(\beF)=0$;
\item[(ii)] 
its restriction to some (general) line $\lda_{gen}\subset\mbb P^3$ is trivializable (hence it is slope semi-stable);
\item[(iii)] 
it satisfies the instanton condition: $H^1(\beF(-2)) = H^2(\beF(-2)) = 0.$
\end{enumerate}\smallskip
We say that an instanton-like vector bundle $\beF$ is \emph{$\cEnd$-general}---it's an \emph{Egen-instanton}, for short\footnote{In a previous version the authors proposed the name \emph{special instantons}. However, this term already appeared in~\cite{beo+fra,hirs+nars}, and the corresponding bundles form a closed subvariety of the full moduli space. In contrast, our condition is an open.}---if, in addition, it satisfies: 
\begin{enumerate}[leftmargin=7ex]
\item[(iv)] 
$H^1\big(\cEnd(\beF)(-2)\big)=H^2\big(\cEnd(\beF)(-2)\big)=0.$
\end{enumerate}
\end{m-definition}

\begin{m-remark}\label{rmk:FonP3}
As mentioned before, the case $r=2$ with the conditions~(i)-(iii), corresponds to the mathematical ('t~Hooft) instanton bundles and it has been intensely investigated. Note that in this case, by the Grauert-M\"ullich theorem, semi-stable vector bundles with $c_1=0$ automatically satisfy the property~(ii).\smallskip 

For arbitrary $r$, Barth-Hulek~\cite{brth+hulk} showed that vector bundles $\beF$ satisfying (i)-(iii) above can be written as the cohomology of a linear monad 
$$
\underbrace{H^2(\beF(-3))}_{\cong\Bbbk^n}\otimes\eO_{\mbb P^3}(-1)
\to\eO_{\mbb P^3}^{\oplus r+2n}\to
\underbrace{H^1(\beF(-1))}_{\cong\Bbbk^n}\otimes\eO_{\mbb P^3}(1).
$$
This construction will be carefully analysed in the next section. 

The condition (iv) is indeed a generic property. By taking the cohomology of the display of the monad, twisted by $\beF^\vee(-2)$, one obtains the exact sequence involving $\beE:=\cEnd(\beF)$: 
$${\scalebox{.85}{$
0\to H^1(\beE(-2))\to H^2(\beF(-3))\otimes H^2(\beF^\vee(-3))\to H^1(\beF(-1))\otimes H^1(\beF^\vee(-1))\to H^2(\beE(-2))\to 0.
$}}$$
Hence the Egen-condition amounts to saying that the middle arrow is an isomorphism; clearly, the relevant vector spaces have the same dimension.
\end{m-remark}

\begin{m-notation}\label{not:moduli} 
We consider the following quasi-projective varieties:
\begin{enumerate}[leftmargin=5ex]
\item 
$\EgI_{\mbb P^3}(r;n),$ the moduli space of Egen-instantons. It is an open subset of the moduli space of slope semi-stable sheaves on $\mbb P^3$.
\smallskip
\item 
$\EgI_{\mbb P^3}(r;n)_\lda,$ the open subspace corresponding to bundles which are trivializable along the line $\lda\subset\mbb P^3$.
\smallskip 
\item 
Let $\bar M_{\mbb P^2}(r;n)$ (resp. $\bar M_{\mbb P^2}(r;n)_{\rm line}$) be the moduli space of rank-$r$, slope semi-stable (resp. framed) vector bundles on $\mbb P^2$, with $c_1=0, c_2=n$. They are $(2rn-r^2+1)$- (resp. $2rn$)-dimensional, cf.~\cite[\S3]{hlc-sheaves}. For simplicity, we call such vector bundles \emph{$\mbb P^2$-instantons}.
\smallskip 
\item 
For two $2$-planes $D,H\subset\ppp$ and $\lda:=D\cap H$, we denote 
$$
\bar M_{D\cup H}(r;n)_\lda:=\frac{\bar M_D(r;n)_\lda\times \bar M_H(r;n)_\lda}{\PGL(r)}
$$
the variety of pairs of bundles, modulo the diagonal action on the framings along $\lda$. 
\item[] 
Furthermore, we denote by 
$$\Theta:\EgI_\ppp(r;n)_\lda\to\bar M_{D\cup H}(r;n)_\lda$$ 
the morphism which sends a vector bundle $\beF$ to its restriction $\beF_{D\cup H}:=\beF\otimes\eO_{D\cup H}$. 
\end{enumerate}
\end{m-notation}

\begin{m-proposition}\label{prop:EI}
\begin{enumerate}[leftmargin=5ex]
\item 
The moduli spaces $\EgI_\ppp(r;n)$ are non-empty, for $n\ges r\ges2$. 
\item 
For any $\beF\in\EgI_\ppp(r;n)$, we have $H^2(\cEnd(\beF))=0$, so its deformations are unobstructed. 
\item 
The differential of the morphism $\Theta$ is an isomorphism everywhere, so $\Theta$ is an \'etale map.  
\item 
Each irreducible component of $\EgI_\ppp(r;n)$ has the expected dimension and the locus corresponding to stable bundles is dense.
\end{enumerate}
\end{m-proposition}

\begin{m-proof}
(i) For $r=2$, consider the union $Z$ of $n+1$ disjoint lines in $\ppp$. The rank-$2$ vector bundle given by the Hartshorne-Serre construction along $Z$ fits into the exact sequence 
$$
0\to\eO_\pi(-1)\to\beF_2\to\cal I_Z(1)\to 0;\quad\beF_2\otimes\eO_Z\cong\eO_Z^{\oplus 2},
$$ 
and one easily verifies that $\beF_2$ is an Egen-instanton. This construction already appears in~\cite[Example 3.1.1]{hart-vb-P3}

For $r>2$, we observe that the rank-$r$ vector bundle $\euf F_r:=\euf O_Y^{\oplus (r-2)}\oplus\,\euf F_2$ still satisfies the conditions~\ref{def:FonP3}. Its general deformation is stable, by (iv) below.

\nit(ii) Take the long exact sequences in cohomology determined by 
$$
\eO_\ppp(-2)\hra\eO_\ppp(-1)\surj\eO_H(-1),\qquad \eO_\ppp(-1)\hra\eO_\ppp\surj\eO_H
$$
twisted by $\cEnd(\beF)$ and use the semi-stabilty of $\beF_H$. 

\nit(iii) The differential of $\Theta$ at $\beF$ is the homomorphism $H^1(\cEnd(\beF))\to H^1(\cEnd(\beF_{D\cup H}))$. The Egen-condition shows that it is indeed an isomorphism. 

\nit(iv) Since $\Theta$ is \'etale, its restriction to each component of $\EgI_\ppp(r;n)$ is dominant. But the stable vector bundles are dense in $\bar M_{D\cup H}(r;n)$ (cf.~\cite[Theorem 3.6]{hlc-sheaves}) and $\beF$ on $\ppp$ is stable as soon as its restriction to $D\cup H$ is so. This shows that the stable bundles $\beF\in\EgI_\ppp(r;n)$ are dense; at such a point, the moduli space is smooth and has the expected dimension.
\end{m-proof}

\begin{m-remark}\label{rmk:lda}
For $r\ges2$, the examples $\beF_r$ above possess the following property: through any point $y\in\ppp$ passes a trivialising line $\lda$ that is, $\beF_r\otimes\eO_\lda\cong\eO_\lda^{\oplus r}$. Thus the same property holds for the general vector bundle in the irreducible component of $\EgI_\ppp(r;n)$ containing this point. Below we show that this holds generally. 

In the next section, such trivialising lines will be used to determine the structure of certain homomorphisms which appear in the monad construction. 
\end{m-remark}

\begin{m-proposition}
Let $\beF$ be an instanton-like vector bundle. Then, for any point $p\in\ppp$, the general line $\lda$ passing through it trivializes $\beF$.
\end{m-proposition}

\begin{m-proof}
Consider the blow-up $\si:Y\to\ppp$ at $p$; it admits a natural projection $\pi:Y\to\mbb P^2$. Let $E=\si^*\eO_\ppp(1)\otimes\pi^*\eO_{\mbb P^2}(-1)$ be the exceptional divisor. 

\nit\underbar{\textit{Claim}}\quad $\si^*\beF$ is slope semi-stable for the paring with 
$$
\pi^*\eO_{\mbb P^2}(1)\cdot(\si^*\eO_\ppp(1)+c\pi^*\eO_{\mbb P^2}(1)),
$$ 
for any $c>0$. (See~\cite{hlc-sheaves} for the terminology.) 

Indeed, let $\tld\beG\subset\si^*\beF$ be a saturated subsheaf, so $\si^*\beF/\tld\beG$ is torsion-free; thus $\tld\beG$ is reflexive and its singularities (if any) are punctual. Then $\det\tld\beG=\si^*\eO_\ppp(k)\otimes\eO_Y(lE)$, for some integers $k,l$, and $\det\tld\beG\subset\si^*\big(\ouset{}{g}{\bigwedge}\beF\big)$; here $g$ stands for the rank of $\tld\beG$. 

Note that  $\si^*\beF/\tld\beG\otimes\eO_E$ is a quotient of $\si^*\beF\otimes\eO_E\cong\eO_E^{\oplus r}$. But $\eO_E(-E)$ is ample on $E$ and $E^3=1$; the semi-stability of $\eO_E^{\oplus r}$ yields $-l\ges0$.  Now we show that $k\les0$, too. Indeed, $\si_*\det\tld\beG=\eO_\ppp(k)\otimes\eI_p^{-l}\subset\ouset{}{g}{\bigwedge}\beF$ and the latter is semi-stable. (We denoted by $\eI_p$ the sheaf of ideals of $p$.) 
It remains to compute
$$
\begin{array}{l} 
c_1(\tld\beG)\cdot\pi^*\eO_{\mbb P^2}(1)\cdot(\si^*\eO_\ppp(1)+c\pi^*\eO_{\mbb P^2}(1)) 
= k\cdot(1+c)+lc\les0.
\end{array}
$$
so $\si^*\beF$ is indeed semi-stable. Finally, we apply~\cite[Theorem A]{hlc-sheaves} and deduce that the restriction of $\si^*\beF$ to the general fibre of $\pi$---corresponding to the general line passing through $p$---is semi-stable (with $c_1=0$), hence it is trivializable. 
\end{m-proof}

%%%%%%%%%%%%%%%%%%%%%%%%%
%%%%%%%%%%%%%%%%%%%%%%%%%

\section{Computations with monads}\label{sct:=0}

Throughout this section, for shorthand, the symbol `$\hra$' indicates a monomorphism and `$\surj$' an epimorphism;  we denote short exact sequences by $A\hra C\surj B$. 

\subsection{The monad construction} 
Let $\beF,\beG$ be instanton-like vector bundles as in Definition~\ref{def:FonP3}. We recall~\cite[p.~340]{brth+hulk} that $\beF^\vee$ is the cohomology of a monad with display 
$$
\scalebox{.99}{$
\xymatrix@R=1.5em@C=3em{
V\otimes\eO(-1)\ar@{^(->}[r]^-{\veps_\beF}\ar@{=}[d]&
\;\eK_\beF\ar@{->>}[r]\ar@{^(->}[d]^{\ev_K^\vee}&
\;\beF^\vee\ar@{^(->}[d]&\qquad(\varstar)
\\ 
V\otimes\eO(-1)\ar@{^(->}[r]^-{\wtld\veps_\beF}&
C\otimes\eO\ar@{->>}[r]^-{\ev_Q}\ar@{->>}[d]^-{\wtld q_\beF}&
\eQ_\beF\ar@{->>}[d]^-{q_\beF}&
\\ 
&W\otimes\eO(1)\ar@{=}[r]&W\otimes\eO(1).&
}
$}
$$ 
If $\beF$ is stable, the monad is uniquely defined up to a natural $\Gl(V)\times\Gl(C)\times\Gl(W)$-action; $V,W,C$ are vector spaces of dimensions $n,n,r+2n$, respectively. In what follows, we enumerate several of its properties. 

\begin{enumerate}[leftmargin=5ex]
\item[(i)] 
There are \emph{canonical} identifications: 

\item[] 
$V\cong V_\beF:=H^2(\beF^\vee(-3)),\;W\cong W_\beF:=H^1(\beF^\vee(-1));$ $H^0(\eQ_\beF) \cong C \cong H^0(\eK_\beF^\vee)^\vee.$

\item[] 
So $\ev_Q$ is $H^0(\eQ_\beF)\otimes\eO_\ppp\to\eQ_\beF$ and $\ev_K^\vee$ is the dual of the evaluation map of $\eK^\vee$.

\item[(ii)] 
The diagram is obtained as follows:
   \begin{itemize}
      \item[(a)] 
      The middle column is the extension defined by $\bone_{W}\in\End(W)=\Ext^1(W(1),\eK_\beF).$ 
      \item[(b)] 
      The middle row is defined by $\bone_{V}\in\End(V)=\Ext^1(\eQ_\beF,V(-1)).$
   \end{itemize}
Although superfluous, we recall~\cite[p. 339]{brth+hulk} the following detail in order to justify (at least morally) the appearance of the Koszul-type isomorphisms in the sequel. The Koszul resolution of a line in $\ppp$ is used to show that the module $\ouset{l\ges 0}{}{\bigoplus}H^1(\beF^\vee(l-1))$ is generated by $H^1(\beF^\vee(-1))$ over the ring $\ouset{l\ges 0}{}{\bigoplus}H^0(\eO_\ppp(l))$. Therefore the top line and the rightmost column of $(\varstar)$, as defined by (a), (b) above, coincide with the minimal resolutions given by~[Proposition~1, p.~327], which is at the core of the construction.
\item[(iii)] 
Let $\lda\subset\mbb P^3$ be a trivialising line for $\beF$ and $\eO(-2)\hra\eO(-1)^{\oplus 2}\surj\eI_\lda\,$ be the Koszul resolution of its sheaf of ideals. By applying $\Hom(\cdot,\beF^\vee(-3))$, we obtain 
\\[.5ex]$\null\hfill 
0\to\underbrace{\Ext^1(\eO(-2),\beF^\vee(-3))}_{=H^1(\beF^\vee(-1))=W_\beF}
\;\ouset{\cong}{\kz_{\beF}}{\lar}\;
\underbrace{\Ext^2(\eI_\lda,\beF^\vee(-3))}_{=H^2(\beF^\vee(-3))=V_\beF}\to 0.
\hfill$\\[.5ex] 
On the right-hand side, we used that $\beF_\lda\cong\eO_\lda^{\oplus r}$. The isomorphism $\kz_\beF$ is given by the Yoneda-product with the element of $\Ext^1(\eI_\lda,\eO(-2))$ corresponding to the Koszul resolution, twisted by the identity of $\beF$.

\item[] 
Note that two trivialising lines $\lda,\lda'$ for $\beF$ induce the same map in cohomology:
   \begin{itemize}
      \item if $\lda,\lda'$ are disjoint, then use $\eI_{\lda\cup\lda'}\hra\eI_\lda\surj\eO_{\lda'};$
      \item if $\lda\cap\lda'=\{y\}$, then use $\eI_{\lda\cup\lda'}\hra\eI_{\lda'}\surj\eO_\lda(-y).$
   \end{itemize} 
\end{enumerate}

As we mentioned earlier, the monad $(\varstar)$ is defined up to a group action; in order to perform computations, we need explicit representatives for its entries. Let $\pl,\pr:\ppp\times\ppp\to\ppp$ be the projections onto the first and second factors, respectively,  and denote by $\Delta\subset\ppp\times\ppp$ the diagonal. The natural evaluation map 
$$
\ev_\ppp:\pl_*(\eI_\Delta\otimes\pr^*\eO(1))\boxtimes\eO(-1)\to \eI_\Delta
$$
is surjective; in fact, 
$$
\pl_*(\eI_\Delta\otimes\pr^*\eO(1))\cong\Omega^1_\ppp(1)
$$ 
and $\ev_\ppp$ is the first term of Beilinson's resolution of the diagonal in $\ppp\times\ppp$. Let us denote 
$$
\tilde\eS:=\Ker(\ev_\ppp)\cong \frac{\Omega^2_\ppp(2)\boxtimes\eO_\ppp(-2)}{\eO_\ppp(-1)\boxtimes\eO_\ppp(-3)}
\cong\frac{\eT_\ppp(-2)\boxtimes\eO_\ppp(-2)}{\eO_\ppp(-1)\boxtimes\eO_\ppp(-3)},
$$
so we get the resolution
$\;\tilde\eS\hra \Omega_\ppp^1(1) \boxtimes\eO_\ppp(-1)\surj\eI_\Delta.$

\begin{m-lemma}\label{lm:hom}
Suppose $\eF$ is stable. Then the display $(\varstar)$ is obtained by applying suitable $\cHom(\cdot,\cdot)$ functors to the diagram below, which is \underbar{\emph{independent}} of the vector bundle $\beF$: 
\begin{m-eqn}{
\scalebox{.99}{$
\xymatrix@R=1.5em@C=1.35em{
\pr^*\eO_\ppp(1)\ar@/^3ex/@{->>}[rr] &\,\eI_\Delta\otimes\pr^*\eO_\ppp(1)\ar@{_(->}[l]&\eO_\ppp(1) 
&\genfrac{}{}{0pt}{1}{\text{apply}}{\cHom(\,\cdot\,,\pl^*\beF^\vee(-3))} 
\\ 
\tilde\eS\otimes\pr^*\eO_\ppp(1)\,\ar@{^(->}[r]
&\pl^*\Omega^1_\ppp(1)\ar@{->>}[r]\ar@{->>}[u]
&\eI_\Delta\otimes\pr^*\eO_\ppp(1)\ar@{^(->}[d] 
&\genfrac{}{}{0pt}{1}{\text{apply}}{\cHom(\pl^*\beF(1),\,\cdot\,)}
\\ 
&\tilde\eS\otimes\pr^*\eO_\ppp(1)\ar@{^(->}[u]&\pr^*\eO_\ppp(1)\ar@/_3ex/@{->>}[uu] 
& 
\\ 
&\genfrac{}{}{0pt}{1}{\text{apply}}{\cHom(\,\cdot\,,\pl^*\beF^\vee(-3))} 
&\genfrac{}{}{0pt}{1}{\text{apply}}{\cHom(\pl^*\beF(1),\,\cdot\,)}
&
}$}
}\label{eq:diag}
\end{m-eqn}
Pointwise, over $y\in\ppp$, one has to restrict the diagram to $\ppp\times\{y\}$. As a consequence, the various homomorphism in $(\varstar)$ are actually induced from the diagram above. 
\end{m-lemma}

\begin{m-proof}
By applying the indicated functors, we obtain 
\begin{m-eqn}{
\scalebox{.77}{$
\xymatrix@R=2.5em@C=1em{
H^2(\beF^\vee(-3))\otimes\eO_\ppp(-1)\ar@{^(->}[r]
&
\mbox{$\underbrace{
\cExt^2_{\pr}\big(\eI_\Delta,\pl^*\beF^\vee(-3)\big)\otimes\eO(-1)}_{=:\;\euf A}
$}
\ar@{->>}[r]\ar@{^(->}[d]^-{a^\vee}
&
\cExt^3_\pr(\eO_\Delta,\beF^\vee(-4)){=}\beF^\vee\kern-1.1ex=\pr_*\pl^*\beF^\vee
\ar@{^(->}[d]
\\ 
R^1\pr_*(\pl^*\eF^\vee(-1)\,{\otimes}\,\tld\eS\,{\otimes}\,\pr^*\eO_\ppp(1) )
\ar@{^(->}[r]
&
\genfrac{}{}{0pt}{1}
{H^1(\ppp,\,\beF^\vee(-1)\otimes\Omega_\ppp^1(1)\,)\,\otimes\,\eO_\ppp}
{\Ext^2_\ppp(\Omega_\ppp^1(1),\beF^\vee(-3))\,\otimes\,\eO_\ppp}
\ar@{->>}[r]^-{b}\ar@{->>}[d]
&
\mbox{$\underbrace{
R^1\pr_*\big(\eI_\Delta\otimes\pl^*\beF^\vee(-1)\big)\otimes\eO(1)}_{=:\;\euf B}
$}
\ar@{->>}[d]
\\ 
&
\cExt^2_\pr(\tld\eS\otimes\pr^*\eO_\ppp(1),\pl^*\eF^\vee(-3) )
&
H^1(\beF^\vee(-1))\otimes\eO_\ppp(1)
}
$}
}\label{eq:beil}
\end{m-eqn}
where $\cExt_\pr$ stands for the relative $\Ext$-functor. 

The top and rightmost extensions are clearly defined by the identity elements in $\End(V_\beF)$ and $\End(W_\beF)$, so $\euf A, \euf B$ coincide with $\eK_\beF$ and $\eQ_\beF$, respectively. The fact that $a,b$ are indeed evaluations morphisms follow from the identities: 
$$
\begin{array}{rl}
H^0(\euf B)
&=
H^0(\ppp,\,R^1\pr_*(\eI_\Delta\otimes\pl^*\beF^\vee(-1))\otimes\eO(1)\,)
\\[1ex] 
=H^1(\ppp\times\ppp,\,\eI_\Delta\otimes(\beF^\vee(-1)\boxtimes\eO(1))\,)
&
=H^1(\ppp,\,\beF^\vee(-1)\otimes\pl_*(\eI_\Delta\otimes\pr^*\eO(1))\,)
\\[1ex] 
&=H^1(\ppp,\,\beF^\vee(-1)\otimes\Omega_\ppp^1(1));
\\[2ex]
H^0(\euf A^\vee)^\vee=\ldots&= \Ext^2_\ppp(\pl_*(\eI_\Delta\otimes\pr^*\eO(1)),\beF^\vee(-3))
\\[1ex]
&=\Ext^2_\ppp(\Omega_\ppp^1(1),\beF^\vee(-3)).
\end{array}
$$
One should also verify that the central term $C$ in $(\varstar)$, coincide with the central terms in the previous diagram. For this, note that we have
$$
0\to C=H^0(\eQ_\beF)\to W_\beF\otimes H^0(\eO_\ppp(1))\to H^1(\eF^\vee),
$$
so $C$ is the kernel of the (pairing) homomorphism on the right. But this pairing is obtained by tensoring $\eI_\Delta\hra\eO_{\ppp\times\ppp}\surj\eO_\Delta$ with $\eF^\vee(-1)\boxtimes\eO(1)$. It follows that $C$ in $(\varstar)$ can indeed be naturally identified with $H^1(\ppp,\,\beF^\vee(-1)\otimes\Omega_\ppp^1(1)\,)$. For $\eK_\beF^\vee$, the argument is similar. 

Now we check, respectively, the isomorphism between the two entries in the leftmost and the  bottom row. The Euler sequence on $\ppp$ and the exact sequence 
$$
\eO_\ppp(-1)\boxtimes\eO_\ppp(-3)\hra \eT_\ppp(-2)\boxtimes\eO_\ppp(-2)\surj \tld\eS
$$ 
yield 
$$\scalebox{.85}{$
\begin{array}{l}
R^1\pr_*(\pl^*\eF^\vee(-1)\otimes\tld\eS\otimes\pr^*\eO_\ppp(1) )\cong H^1(\eF^\vee{\otimes}\eT_\ppp(-3))\otimes\eO_\ppp(-1)\cong H^2(\eF^\vee(-3))\otimes\eO_\ppp(-1).
\\[1.5ex] 
\cExt^2_\pr(\tld S\otimes\pr^*\eO_\ppp(1),\pl^*\beF^\vee(-3))\cong H^2(\Omega^1_\ppp(1)\otimes\beF^\vee(-2))\otimes\eO_\ppp(1)\cong H^1(\beF^\vee(-1))\otimes\eO_\ppp(1).
\end{array}
$}$$
The verifications done so far already ensure that our diagram agrees with the display $(\varstar)$; indeed, the latter is determined by expanding either the rightmost column, through the isomorphism  $\Ext^1(W(1),\eK_\beF)\srel{\cong}{\to}\Ext^1(W(1),\beF^\vee)$, or the top row, through $\Ext^1(\eQ_\beF,V(-1))\srel{\cong}{\to}\Ext^1(\beF^\vee,V(-1))$. However, it is worth clarifying the remaining (rather mysterious) coincidence of the central terms. This can be seen again by restricting to a hyperplane $H\subset\ppp$. The Euler sequence implies: 
$$
\begin{array}{lcl}
H^1(\beF^\vee(-2)\otimes\Omega^1_\ppp(1))=0,&& H^1(\beF^\vee(-1))\srel{\cong}{\to} H^2(\beF^\vee(-2)\otimes\Omega^1_\ppp(1));
\\[1.5ex]
H^1(\beF^\vee(-2)\otimes\eT_\ppp(-1))\srel{\cong}{\to}H^2(\beF^\vee(-3)),&&H^2(\beF^\vee(-2)\otimes\eT_\ppp(-1))=0.
\end{array}
$$
The commutative diagrams below (tensored by $\beF^\vee(-1)$ and $\beF^\vee(-2)$, respectively)\\[1ex]  
\scalebox{.85}{$\xymatrix@R=1em@C=1em{
\Omega^1_\ppp\ar@{^(->}[r]\ar[dd]&\Omega^1_\ppp(1)\ar@{->>}[r]\ar[d]&\Omega^1_\ppp(1)_H=\Omega^1_H(1)\oplus(\eI_H/\eI_H^2)(H) \ar[dd]
\\ 
&H^0(\eO_\ppp(1))\otimes\eO_\ppp\ar[d]&
\\ 
\eO_\ppp(-H)\ar@{^(->}[r]&\eO_\ppp\ar@{->>}[r]&\eO_H 
}\hfill$}\\[1ex]
\null\hfill\scalebox{.85}{$\xymatrix@R=1em@C=1em{
\eO_\ppp(-1)\ar@{^(->}[r]\ar@{=}[dd]&H^0(\eO_\ppp(1))^\vee\otimes\eO_\ppp\ar@{->>}[r]\ar[dd]&\eT_\ppp(-1)\ar@{->>}[d]
\\ 
&&\eT_\ppp(-1)_H=\eT_H(-1)\oplus\eO_H\ar[d]
\\
\eO_\ppp(-1)\ar@{^(->}[r]&\eO_\ppp\ar@{->>}[r]&\eO_H
}$}\\[1ex	] 
yield the exact sequences 
$$\scalebox{.85}{$
\begin{array}{c}
0\to H^1(\beF^\vee(-1)\otimes\Omega^1_\ppp(1))\to H^1(\beF^\vee_H\otimes\Omega^1_H)\oplus 
\underbrace{H^1(\beF^\vee_H(-1))\to H^1(\beF^\vee(-1))}_\cong,
\\[5ex]
\underbrace{H^2(\beF^\vee(-3))\to H^1(\beF^\vee_H(-2))}_\cong\oplus H^1(\beF^\vee_H\otimes\eT_H(-3))\to H^2(\beF^\vee(-3)\otimes\eT_\ppp(-1))\to0,
\end{array}
$}$$
whose rightmost (resp. leftmost) terms are isomorphic. Since $\eT_H(-1)\cong\Omega^1_H(2)$, it follows that $H^1(\beF^\vee(-1)\otimes\Omega^1_\ppp(1))$ and $H^2(\beF^\vee(-3)\otimes\eT_\ppp(-1))$ are isomorphic and both can be identified with $H^1(\beF^\vee_H\otimes\Omega^1_H)$.
\end{m-proof}

\begin{m-lemma}\label{lm:loc-inv}
Let $\beF$ be an instanton-like vector bundle. Take a point $y\in\ppp$ and consider a line $\lda\subset\ppp$ passing through it, such that $\beF_\lda\cong\eO_\lda^{\oplus r}$.  Then the inclusion $\eI_\lda\subset\eI_y$ determines a (local) inverse to $\veps_\beF$ that is, the choice of $\lda$ determines a splitting of the top row of $(\varstar)$ at $y$.
\end{m-lemma}

\begin{m-proof}
By applying the functor $\Hom(\cdot,\beF^\vee(-3))$ to $\eI_\lda\hra\eI_y\surj\eO_\lda(-y)$, we obtain 
$$
0\to\underbrace{\Ext^2(\eO_\lda(-y),\beF^\vee(-3))}_{=\Ext^3(\eO_y,\beF^\vee(-3))=\beF^\vee(1)_y}
\to 
\underbrace{\Ext^2(\eI_y,\beF^\vee(-3))}_{=\eK_\beF(1)_y}
\to 
\underbrace{\Ext^2(\eI_\lda,\beF^\vee(-3))}_{=H^2(\beF^\vee(-3))=V_\beF}\to0.
$$
This yields a (local) inverse of $\veps_\beF$, due to the commutativity of 
$$
\xymatrix@R=1.5em{\eI_\lda\ar[r]\ar[d]&\eI_y\ar[d]\\ \eO_\ppp\ar@{=}[r]&\eO_\ppp.}
$$
Clearly, the construction can be done locally about $y$ and we obtain local splittings (over open subsets $\cU\subset\ppp$) denoted by $\eta_\cU:\eK_{\beF\;\rst\cU}\to V_\beF(-1)_\cU$. We can assume that they are small enough so that the middle column splits over them, too.
\end{m-proof}

%%%%%%%%%%%%%%%%%%%%%%%%%%%%%%%%
%%%%%%%%%%%%%%%%%%%%%%%%%%%%%%%%

\subsection{Computations} 

Henceforth it is convenient to bear in mind that cohomology classes can be represented as {\v C}ech cocycles, which are genuine sections over open subsets, so one can understand easier the effect of homomorphisms on them. It is common in the literature to denote cocycles by $\cal Z^\bullet(\,\cdot\,)$. 

Our next goal is to prove the commutativity of the diagram:
$$\scalebox{.9}
{$\xymatrix@R=1.5em@C=1.35em{
&&&
\null\hspace{3em} W_\beF\otimes H^1(\beG(-1))\hspace{3em} \null 
\ar[d]^-{\beta_\beF\otimes\bone_\beG}_-\cong
\ar[dl]_-{\kz_\beF\otimes\kz_\beG}^(.37)\cong 
\\ 
\mbox{$\underbrace{H^1(\eK_\beF\otimes\beG(-2))}_{=0}$}\ar[r] 
& 
H^1(\cHom(\beF,\beG)(-2))\ar[r] 
& 
V_\beF\otimes H^2(\beG(-3))\ar[r]_-{H^2(\veps_\beF\otimes\bone_\beG)}&
H^2(\eK_\beF\otimes\beG(-2))
}$
}$$

\begin{m-proposition}\label{prop:=}
For general Egen-instantons, the diagram above is commutative that is, 
\begin{m-eqn}{
H^2(\veps_\beF\otimes\bone_\beG)\circ(\kz_\beF\otimes\kz_\beG)
=\beta_\beF\otimes\bone_\beG.
}\label{eq:comm}
\end{m-eqn}
Therefore, it holds $H^1(\cHom(\beF,\beG)(-2))=H^2(\cHom(\beF,\beG)(-2))=0$.
\end{m-proposition}

Note that the homomorphisms on both sides of the equation above act separately on $\beF,\beG$. If the rank is allowed to vary, one can say that the left- and the right-hand side are bilinear functors (from the category of semi-stable vector bundles, with $c_1=c_3=0$, to the category of homomorphisms between vector spaces). One may turn the statement around:  $(\beta_\beF\otimes\bone_\beG)^{-1}\circ H^2(\veps_\beF\otimes\bone_\beG):V_\beF\otimes V_\beG\to W_\beF\otimes W_\beG$ is a (functorial and bilinear) isomorphism, so it's natural to ask what is this map. 

This leads to the idea to analyse the effect on $\beF$ and $\beG$ separately. The $\beG$-component is easier to understand---it is obtained by tensoring with the identity map of $\beG$---, while the $\beF$-component is rooted deeper into the structure of the monad. 

\begin{m-lemma}\label{lm:commG}
Let $\beF,\beG$ be two stable instanton-like vector bundles. Then it holds:  
\begin{m-eqn}{
(\beta_\beF\otimes\bone_\beG)^{-1}\circ H^2(\veps_\beF\otimes\bone_\beG)=\chi_\beF\otimes\kz_\beG^{-1},
\quad\text{with}\; \chi_\beF\in\Hom(\,H^2(\beF(-3)),H^1(\beF(-1))\,). 
}\label{eq:G}
\end{m-eqn}
That is, the left-hand side acts on the $\beG$-component of the tensor product the same way as the inverse of the Koszul map. 
Moreover, the homomorphism $\chi_\beF$ depends only of $\beF$ (it is independent of $\beG$).
\end{m-lemma}
Note that we don't require $\beF$ or $\beG$ to be Egen-instantons. Also, in the proof, we shall use no information about the structure of $V_\beF,W_\beF$, but only the properties of $\beG$.

In order to clarify our reasoning, let us point out two rearrangements in the display $(\varstar)$, which should clarify our approach involving restrictions to $2$-planes in $\ppp$:
$$
V_\beF\otimes\beG(-3)=V_\beF(-1)\otimes\beG(-2)\quad\text{and}\quad W_\beF(1)\otimes\beG(-2)=W_\beF\otimes\beG(-1).
$$ 
They are necessary to apply $\veps_\beF,\beta_\beF$, respectively, and correspond to division (for $V_\beF,W_\beF$) by a linear equation and multiplication (for $\beG$) by the same factor.

\begin{m-proof} 
Our reasoning involves three steps: we start by analysing the structure of the homomorphisms induced by $\veps_\beF$ and by $\beta_\beF$; finally we compose the two maps. \medskip 

\nit\underbar{\textit{Step 1}}\quad Let $\lda$ be a trivialising line for $\beG$ and $D,H$ two planes containing it. The restrictions $\beG_D,\beG_H$ are automatically semi-stable and, for generic choices, they are actually stable. For analysing the homomorphism induced by $\veps$, we consider the diagram:
$$\scalebox{.85}{$
\xymatrix@R=1.5em@C=2em{
V_\beF\otimes H^1(\beG(-2)_H)\ar@/_15ex/[ddd]^-\cong\ar@{^(->}[r]
&
H^1(\eK_\beF\otimes\beG(-1)_H)\ar@{^(->}[r]
&
C_\beF\otimes H^1(\beG(-1)_H)\ar@{=}[r]^{\text{restrict}}_-{\text{to $\lda$}}
&
C_\beF\otimes H^1(\beG(-2)_H)\ar@{->>}[dd]^-\cong
\\ 
V_\beF\otimes H^1(\beG(-2)_{D\cup H})\ar@{->>}[u]\ar[d]^-\cong\ar[r]
&
H^1(\eK_\beF\otimes\beG(-1)_{D\cup H})\ar[d]\ar[r]\ar[u]
&
C_\beF\otimes H^1(\beG(-1)_{D\cup H})\ar@{->>}[u]\ar@{->>}[dr]
&
\\ 
V_\beF\otimes H^2(\beG(-4))\ar@{->>}[d]\ar[r]
&
H^2(\eK_\beF(1)\otimes\beG(-4))\ar@{->>}[d]\ar[rr]^-\cong
&&
C_\beF\otimes H^2(\beG(-3))
\\ 
V_\beF\otimes H^2(\beG(-3))\ar[r]_-{H^2(\veps_\beF\otimes\bone_\beG)}\ar@{..>}[rrru]|(.6){\;\delta\;}
&
H^2(\eK_\beF(1)\otimes\beG(-3))
&
}
$}$$
We observe that it is commutative: this fact can be seen by moving along the second row. 

First we claim that the dotted arrow $\delta$---defined by following the top row---acts as the identity on $H^2(\beG(-3))$. Indeed, a lifting of an element in $H^2(\beG(-3))$ to $H^2(\beG(-4))$ amounts to dividing the corresponding cocycle by a linear equation (the surjectivity of the arrow ensures that such a division makes sense). The homomorphism $\veps_\beF$ is a linear combination of elements in $H^0(\eO_\ppp(1))$, so $\veps_\beF=\ouset{j=0}{3}{\sum} c_j L_j$, where $L_0,\dots,L_3$ is a basis of linear forms and $c_j\in V_\beF^\vee\otimes C_\beF$. Then, by following the third row, we deduce that $\delta$ acts on $v\otimes h\in V_\beF\otimes H^2(\beG(-3))$ as follows: for $j=0,\dots,3$, there is a representative $\tld h_j\in\cal Z^2(\beG(-3))$ of $h$, such that the quotient $\frac{\tld h_j}{L_j}\in\cal Z^2(\beG(-4))$ is well-defined, so we have 
$$
\delta(v\otimes h)=\text{the cohomology class defined by}\;\;\ouset{j=0}{3}{\sum} c_j(v)\otimes L_j\cdot \frac{\tld h_j}{L_j}\;=\Big(\ouset{j=0}{3}{\sum} c_j(v)\Big)\otimes h.
$$
The linear factor required for lifting $h$ to $H^2(\beG(-4))$ cancels out by applying $\veps_\beF$, so the only operation performed on the $H^2(\beG(-3))$-factor is the tensor product by elements of $C_\beF$. It is also apparent that the $\beF$-component of $\delta$ is determined only by $\beF$. \medskip 

\nit\underbar{\textit{Step 2}}\quad Now we turn our attention to $\beta_\beF$: by tensoring the exact sequences 
$$
\eK_\beF(-1)\hra C_\beF(-1)\surj W_\beF\otimes\eO_\ppp
\quad\text{and}\quad 
\beG(-2)\hra\beG(-1)\surj\beG(-1)_H
$$ 
we obtain the diagram 
$$\scalebox{.85}{$
\xymatrix@R=1.5em@C=1.35em{
\eK_\beF\otimes\beG(-3)\ar@{^(->}[r]\ar@{^(->}[d]
&
\eK_\beF\otimes\beG(-2)\ar@{^(->}[d]\ar@{->>}[r]
&
\eK_\beF\otimes\beG(-2)_H\ar@{^(->}[d]\ar@{-->}[d]
\\ 
C_\beF\otimes\beG(-3)\ar@{^(->}[r]\ar@{->>}[d]
&
C_\beF\otimes\beG(-2)\ar@{->>}[d]\ar@{->>}[r]\ar@{..>}[dr]|{\quad m\quad}
&
C_\beF\otimes\beG(-2)_H\ar@{->>}[d]
\\ 
W_\beF\otimes\beG(-2)\ar@{^(->}[r]
&
W_\beF\otimes\beG(-1)\ar@{->>}[r]
&
W_\beF\otimes\beG(-1)_H.
}
$}$$
Then $\euf M:=\Ker(m)$ satisfies the isomorphism 
$$W_\beF\otimes H^1(\beG(-1)_H)\srel{\cong}{\to} H^2(\euf M)$$ 
and it also fits into the following three diagrams:\\[2ex] 
$
\begin{array}{c|c|c}
\kern-1ex\scalebox{.61}{$\xymatrix@R=1.35em@C=1.35em{
\eK_\beF\otimes\beG(-3)\ar@{^(->}[r]\ar@{^(->}[d]
&
\eK_\beF\otimes\beG(-2)\ar@{^(->}[d]\ar@{->>}[r]
&
\eK_\beF\otimes\beG(-2)_H\ar@{=}[d]
\\ 
C_\beF\otimes\beG(-3)\ar@{^(->}[r]\ar@{->>}[d]
&
\euf M\ar@{->>}[r]\ar@{->>}[d]
&
\eK_\beF\otimes\beG(-2)_H
\\ 
W_\beF\otimes\beG(-2)\ar@{=}[r]
&
W_\beF\otimes\beG(-2)
&
}$}\kern-1ex
& %ENDS THE FIRST ENTRY OF ARRAY
\kern-1ex\scalebox{.61}{$\xymatrix@R=1.35em@C=1.35em{
&
&
\eK_\beF\otimes\beG(-2)_H\ar@{^(->}[d]
\\ 
C_\beF\otimes\beG(-3)\ar@{^(->}[r]\ar@{^(->}[d]
&
C_\beF\otimes\beG(-2)\ar@{=}[d]\ar@{->>}[r]
&
C_\beF\otimes\beG(-2)_H\ar@{->>}[d]
\\ 
\euf M\ar@{^(->}[r]\ar@{->>}[d]
&
C_\beF\otimes\beG(-2)\ar@{->>}[r]
&
W_\beF\otimes\beG(-1)_H
\\ 
\eK_\beF\otimes\beG(-2)_H
}$}\kern-1ex
& %ENDS THE SECOND ENTRY OF ARRAY
\kern-1ex\scalebox{.61}{$\xymatrix@R=1.35em@C=1.35em{
&
&
W_\beF\otimes\beG(-2)\ar@{^(->}[d]
\\ 
\eK_\beF\otimes\beG(-2)\ar@{^(->}[r]\ar@{^(->}[d]
&
C_\beF\otimes\beG(-2)\ar@{=}[d]\ar@{->>}[r]
&
W_\beF\otimes\beG(-1)\ar@{->>}[d]
\\ 
\euf M\ar@{^(->}[r]\ar@{->>}[d]
&
C_\beF\otimes\beG(-2)\ar@{->>}[r]
&
W_\beF\otimes\beG(-1)_H
\\ 
W_\beF\otimes\beG(-2)
}$}% ARRAY, END OF FIRST LINE 
\\ \rm (I)&\rm (II)&\rm (III)
\end{array}
$\\[2ex] 
They all imply the commutativity of the diagram
$$\scalebox{.85}{$
\xymatrix{
H^2(\eK_\beF\otimes\beG(-3))\ar@{->>}[r]\ar@{->>}[d]^-\cong\ar@{}[dr]|{\rm(I)}
&
H^2(\eK_\beF\otimes\beG(-2))\ar@{->>}[d]^-\cong\ar@{}[dr]|{\rm(III)}
&
W_\beF\otimes H^1(\beG(-1))\ar[l]_-\cong\ar[d]^-\cong
\\ 
C_\beF\otimes H^2(\beG(-3))\ar@{->>}[r]\ar@{..>}[urr]|(.35){\;\gamma\;}\ar@{}[dr]|{\rm(II)}
&
H^2(\euf M)
&
W_\beF\otimes H^1(\beG(-1)_H)\ar[l]_-\cong
\\ 
C_\beF\otimes H^1(\beG(-2)_H)\ar@{->>}[u]^-\cong\ar@{->>}[r]
&
W_\beF\otimes H^1(\beG(-1)_H)\ar@{->>}[u]^-\cong\ar@{=}[ur]
&
}
$}$$
By moving along the lower edges of the diagram, we see that the dotted homomorphism $\gamma$ is also a tensor product: its $\beF$-component is a composition of various homomorphisms between cohomology groups---that is, $C_\beF,W_\beF$---determined by $\beF$.  The $\beG$-component acts on $H^2(\beG(-3))$ as the inverse of the Koszul map: indeed, the lower side of the diagram is obtained by tensoring with $\beG(-1)$ and applying $\veps,\veps_C$ to either one of the following: 
$$
\begin{array}{c|c}
\scalebox{.85}{$\xymatrix@R=1.5em@C=1.35em{
\eO_\ppp(-D-H)\ar@{^(->}[d]&&\\ \eO_\ppp(-D)\ar@{->>}[d]&&\\ \eO_H(-\lda)\ar@{^(->}[r]&\eO_H\ar@{->>}[r]&\eO_\lda
}$}
& 
\scalebox{.85}{\xymatrix@R=1.5em@C=1.35em{
&\eO_\ppp(-H)\ar@{=}[r]\ar@{^(->}[d]&\eI_H\ar@{^(->}[d]\\ \eO_\ppp(-D-H)\ar@{^(->}[r]\ar@{=}[d]&\eO_\ppp(-D)\oplus\eO_\ppp(-H)\ar@{->>}[r]\ar@{->>}[d]&\eI_\lda\ar@{->>}[d]\\ \eO_\ppp(-D-H)\ar@{^(->}[r]&\eO_\ppp(-D)\ar@{->>}[r]&\eI_{\lda\subset H}=\eO_H(-\lda).
}}
\end{array}
$$

\nit\underbar{\textit{Step 3}}\quad The composition in~\eqref{eq:G} coincides with $\gamma\circ\delta$, so it is a tensor product of two linear maps, and the $\beG$-component is $\kz_\beG^{-1}$, while the $\beF$-component is independent of $\beG$. 
\end{m-proof}

Henceforth we are interested in determining the homomorphism $\chi_\beF$. Although it is tempting to use the explicit expressions in the previous lemma and in~\eqref{eq:beil}, the multitude of identifications make difficult pursuing this path. Let us recall that a general element $\beF\in\EgI_\ppp(r;n)$ is stable and through each point $y\in\ppp$ passes a line $\lda$ such that $\beF\otimes\eO_\lda\cong\eO_\lda^{\oplus r}$. We are going to consider general elements in this sense. 

\begin{m-lemma}{\label{lm:commF}}
Let $(\varstar')$ be the diagram obtained by replacing in $(\varstar)$ the entry $V_\beF$ by $W_\beF$. (That is, the top and middle horizontal extensions are given by $\kz_\beF$ rather than the identity of $V_\beF$.) Let $\veps'_\beF$ denote the top inclusion, instead of $\veps_\beF$. 

Suppose $\beF,\beG$ are in the same component of $\EgI_\ppp(r;n)$ and general enough. Then, in $(\varstar)$, it holds $\chi_\beF=\kz_\beF^{-1}$. Equivalently, in the diagram $(\varstar')$,  $H^2(\veps_\beF\otimes\bone_\beG)$ is invertible and the $\beF$-component of $H^2(\veps'_\beF\otimes\bone_\beG)^{-1}\circ (\beta_\beF\otimes\bone_\beG)$ is the identity of $W_\beF$.
\end{m-lemma}

\begin{m-proof}
Conceptually, the statement is due to the fact that $(\varstar)$ is obtained by applying $\cHom$-functors to~\eqref{eq:diag} and compositions of arrows correspond to Yoneda-products (concatenations): $\beta_\beF$ in $(\varstar)$ is given by pairing with $e\in\Ext^1(\eI_\Delta,\tilde\eS)$  twisted by the identity of $\beF$. On the other hand, $\kz_\beF$ is the connecting homomorphism for the Koszul resolution. Hence, in order to relate them, one should compare their defining extensions.\medskip 

The Egen-condition implies that $H^1(\cHom(\beF,\beG)(-2))=H^2(\cHom(\beF,\beG)(-2))=0$, for all $(\beF,\beG)$ in some neighbourhood of diagonal in $\EgI_\ppp(r;n)\times\EgI_\ppp(r;n)$. For such pairs, $H^2(\veps_\beF\otimes\bone_\beG)$ is an isomorphism and we can take its inverse. 

Recall that a short exact sequence $A\hra C\surj B$ of locally free sheaves determines a homomorphism $\cal Z^k(B)\to \cal Z^{k+1}(A)$ at the level of cocycles, given by the cup product with a representative in $\cal Z^1(\cHom(B,A))$ of the extension class, which induces the connecting homomorphism in cohomology.\footnote{We are interested in $k=0,1$. For Dolbeault cocycles, our statement is clear: the product of $\dbar$-closed sections is still $\dbar$-closed. For {\v C}ech cocycles, we use standard notation: let $(e_{ij})_{i,j}$ be a representative of the extension class. Then, for $k=0$, define $(b_i)_i\mt (a_{ij}:=e_{ij}(b_{j\rst U_i}))_{i,j}$;\\ for $k=1$, define $(b_{ij})\mt (a_{ijk}:=e_{jk\rst U_i}(b_{jk\rst U_i})+e_{ij\rst U_k}(b_{ij\rst U_k}))_{i,j,k}.$}

Back to our situation, since $H^2(\veps_\beF\otimes\bone_\beG)$ is invertible, each element of $H^2(\eK_\beF\otimes\beG(-2))$ admits a cocycle representative which belongs to the image of $\veps_\beF\otimes\bone_\beG$. Thus the $\beF$-component of $H^2(\veps_\beF\otimes\bone_\beG)^{-1}\circ (\beta_\beF\otimes\bone_\beG)$ is the map $\cal Z^1(\beF^\vee(-1))\to\cal Z^2(\beF^\vee(-3))$, which (locally) is given by the composition of the connecting homomorphism $\beta_\eF$ followed by the projection $\eta_\cU:\eK_{\beF\;\rst\cU}\to V_\beF(-1)_{\cU}$. The projection is induced by the inclusion $\eI_\lda\subset\eI_y$, so this is the same as applying the connecting map of the pull-back (Yoneda pairing) of the extension $\eS\hra\Omega^1(1)\to\eI_y$ by the inclusion $j_\lda:\eI_\lda\hra\eI_y$: 
\begin{m-eqn}{
\scalebox{.85}{$
\xymatrix@R=1.5em@C=1.35em{
\eS\ar@{=}[rr]\ar@{^(->}[d]&&\eS\ar@{^(->}[d]&
\\ 
j_\lda^*\Omega^1_\ppp(1)\ar@{->>}[d]\ar@{^(->}[rr]&&
\Omega^1_\ppp(1)\ar@{->>}[d]\ar@{=}[r]&H^0(\eI_y(1))\otimes\eO_\ppp(-1)
\\ 
\eI_\lda\ar@{^(->}[rr]^-{j_\lda}&&\eI_y&
}
\kern5ex \eS\cong\frac{\Omega^2_\ppp(2)}{\eO_\ppp(-1)}.
$}
}\label{eq:pb}
\end{m-eqn}
We are going to show that the left column is induced by the Koszul resolution of $\eI_\lda$. Note that $h^1(\eS(1))=0$, so we have 
$$
\scalebox{.85}{$
\xymatrix@R=1.5em@C=3em{
H^0(j_\lda^*\Omega^1_\ppp(2))\ar@{^(->}[r]\ar@{->>}[d]&H^0(\Omega^1_\ppp(2))\ar@{->>}[d] \\ 
H^0(\eI_\lda(1))\ar@{^(->}[r]&H^0(\eI_y(1))
}$}
$$
and we obtain the diagram 
\begin{m-eqn}{
\scalebox{.85}{$\xymatrix@R=1.5em@C=3em{
\eO(-2)\,\ar@{^(->}[r]^-{top}\ar@{^(->}[d]&\eS\ar@{^(->}[d]\ar@{->>}[r]
&\displaystyle\frac{\eS}{\eO(-2)}\ar[d]^-\cong
\\ 
\eO(-1)^{\oplus 2}\cong H^0(\eI_\lda(1))\otimes\eO(-1)\,\ar@{^(->}[r]^-{mid}\ar@{->>}[d]
&j_\lda^*\Omega^1_\ppp(1)\ar@{->>}[d]\ar@{->>}[r]&\displaystyle\frac{j_\lda^*\Omega^1_\ppp(1)}{\eO(-1)^{\oplus 2}}
\\ 
\eI_\lda\ar@{=}[r]&\eI_\lda
}$}
}\label{eq:po}
\end{m-eqn}
The map $mid$ is the natural evaluation and $top$ is induced by 
$$\scalebox{.85}{$
H^0(\eI_y(1))\otimes\eO(-1)\to\Omega^1(1)\quad\leadsto\quad 
\ouset{}{2}{\bigwedge}\big(\,H^0(\eI_y(1))\otimes\eO(-1)\,\big)\to\Omega^2(2)\to\eS.
$}$$
This shows that (locally) the push-out by $top$ of the Koszul resolution of $\eI_\lda$ is the left column of~\eqref{eq:pb}. So, at the cohomological level, the globally defined homomorphism $H^2(\veps_\beF\otimes\bone_\beG)^{-1}\circ(\beta_\beF\otimes\bone_\beG)$ acts on $W_\beF$ the same as $\kz_\beF$.  
\end{m-proof}

\begin{m-proof}\hspace{-.75ex}(of Proposition~\ref{prop:=})\quad It remains to combine the previous steps. 

We start by considering the case when $\beF,\beG$ belong to the same irreducible component of $\EgI_\ppp(r;n)$. Suppose moreover that $\beF,\beG$ are sufficiently close to each other, as in Lemma~\ref{lm:commF}, so $H^2(\veps_\beF\otimes\bone_\beG)$ is bijective. Then, for the new diagram $(\varstar')$, we have 
$$
H^2(\veps'_\beF\otimes\bone_\beG)^{-1}\circ(\beta_\beF\otimes\bone_\beG)=\;\bone_\beF\otimes \chi_\beG,
$$ 
for some homomorphism $\chi_\beG$. This $\beG$-component was determined in Lemma~\ref{lm:commG}: it is just the Koszul homomorphism of $\beG$. 

So far we proved that $H^2(\veps_\beF\otimes\bone_\beG)=(\beta_W\otimes\bone_\beG)\circ(\kz_\beF\otimes\kz_\beG)$, for $(\beF,\beG)$ close enough to each other. By keeping $\beF$ fixed and allowing $\beG$ to vary (in the same component as $\beF$), this can be interpreted as the identity of two sections in a $\cHom$-bundle, between the vector bundles whose fibres over $\beG\in\EgI_\ppp(r;n)$ are $V_\beF\otimes H^2(\beG(-3))$ and $H^2(\eK_\beF\otimes\beG(-2))$, respectively. They agree on a neighbourhood of $\beF$, so the identity holds on the irreducible component containing $\beF$, wherever the vector bundles and the homomorphisms make sense. 

Now let $\beF,\beG\in\EgI_\ppp(r;n)$ be arbitrary; the proof of Lemma~\ref{lm:commF} doesn't apply directly. However, we proved that $\chi_\beF$ in~\eqref{eq:G} is independent of $\beG$, it is the composition of various homomorphisms between cohomology groups determined by $\beF$ only. Therefore $\chi_\beF$ can be computed by using an instanton-like vector bundle $\beG'$ in the same component as $\beF$. The conclusion follows from the previous step.
\end{m-proof}

%%%%%%%%%%%%%%%%%%%%%%%%%
%%%%%%%%%%%%%%%%%%%%%%%%%

\section{The main results}\label{sct:main}

Now we are in position to prove the Theorem stated in the Introduction. 

%%%%%%%%%%%%%%%%%%%%%%%%

\subsection{The irreducibility}\label{ssct:irred}

\begin{m-theorem}\label{thm:irred}
$\EgI_\ppp(r;n)$ is irreducible. 
\end{m-theorem}

\begin{m-proof}
Indeed, by Proposition~\ref{prop:EI}, the restriction of the map $\Theta$ to each irreducible component of $\EgI_\ppp(r;n)$ dominates $\bar M_{D\cup H}(r;n)_\lda$. Let $\beF,\beG$ be general Egen-instantons which are mapped to the same (general) point in $\bar M_{D\cup H}(r;n)_\lda$. Then, Proposition~\ref{prop:=} states that $H^1(\cHom(\beF,\beG)(-2))=0$, so the isomorphism between $\beF_{D\cup H}$ and $\beG_{D\cup H}$ lifts to an isomorphism over all $\ppp$. 
\end{m-proof}

This short argument strongly contrasts the lengthy computations~\cite{tikh-odd,tikh-even} in the rank-$2$ case.

%%%%%%%%%%%%%%%%%%%%%%

\subsection{The rationality}\label{ssct:rat}
Let $D,H\subset\mbb P^3$ be a generic pair (wedge) of $2$-planes, intersecting along $\lda:=D\cap H$.

\begin{m-proof} 
We showed in~\ref{prop:EI} that the restriction map $\Theta:\EgI_\ppp(r;n)_\lda\to M_{D\cup H}(r;n)_\lda$ is \'etale. Proposition~\ref{prop:=} implies that an isomorphism between the restrictions to $D\cup H$ of two Egen instantons actually comes from isomorphism over $\ppp$; thus $\Theta$ is birational. It remains to observe that~\cite[Corollary 3.11]{hlc-sheaves} implies that $\bar M_{\mbb P^2}(r;n)_\lda$ is a rational variety. 
\end{m-proof}

We remark that we actually obtain an explicit description of the general Egen-instanton on the projective space.

\begin{m-theorem}\label{thm:B} 
\begin{enumerate}[leftmargin=5ex]
\item 
The assignment $\beF\mt (\beF_D,\beF_H)$ induces the birational map 
$$\EgI_{\mbb P^3}(r;n)_\lda\srel{\Theta}{\lar}\bar M_{D\cup H}(r;n) 
\cong M_{\mbb P^2}(r;n)\times M_{\mbb P^2}(r;n)_{\lda}.$$ 
Hence a general Egen-instanton bundle $\beF$ on $\mbb P^3$ is uniquely determined by its restrictions $(\beF',\beF'')$ to the $2$-planes $D,H$ and the gluing data $\beF'_\lda\cong\eO_\lda^{\oplus r}\cong\beF''_\lda$ (trivializations, up to simultaneous $\PGL(r)$-action).
\item 
Let $Q\cong\mbb P^1\times\mbb P^1\subset\ppp$ be a smooth quadric. The assignment 
$$
\Theta':\EgI_\ppp(r;n)\dashto\bar M_Q(r;2n),\quad\beF\mt\beF_Q
$$ 
yields a birational map to the moduli space of semi-stable bundles on $Q$, with $c_2=2n$.
\end{enumerate}
\end{m-theorem}

\begin{m-proof} 
The first statement is clear. For the second, note that the same argument as in Proposition~\ref{prop:EI} shows that $\Theta'$ is \'etale; by~\ref{prop:=}, it is actually birational. 
\end{m-proof}

\begin{m-remark}
(i) At the infinitesimal, deformation-theoretic level, these facts are reflected in the isomorphism:
$$
H^1(\mbb P^3,\cEnd(\beF))\srel{\cong}{\to} H^1(D\cup P,\cEnd(\beF))
\srel{\cong}{\to} H^1(Q,\cEnd(\beF)).
$$

\nit(ii) The results obtained in~\cite[Theorem 3.6]{hlc-sheaves} yield detailed descriptions of $\bar M_{\mbb P^2}(r;n)$ and $\bar M_Q(r;2n)$: they are irreducible, rational varieties. Their general elements are, respectively, the kernels of surjective homomorphisms: 
$$
\begin{array}{lll}
-\;\text{for}\;\bar M_{\mbb P^2}(r;n),&\quad& 
 {\eI_{p}^a(a)}^{\oplus r-\rho}\oplus{\eI_{p}^{a+1}(a+1)}^{\oplus \rho} \lar\ouset{j=1}{n}{\bigoplus}\;\eO_{l_j}(1),
\\ && 
n=ar+\rho,\;0\les\rho<r,\;\text{and $l_1,\dots,l_n\subset\mbb P^2$}  
\\ && 
\text{is a bouquet of distinct lines passing through the point $p\in\mbb P^2$;} 
\\[2ex] 
-\;\text{for}\;\bar M_{Q}(r;2n),&&
{\eO_{\mbb P^1_{left}}(a')}^{\oplus r-\rho'}\oplus{\eO_{\mbb P^1_{left}}(a'+1)}^{\oplus\rho'} \lar\ouset{j=1}{2n}{\bigoplus}\;\eO_{{x_j}\times\mbb P^1_{right}}(1),
\\ 
Q=\mbb P^1_{left}\times\mbb P^1_{right}&\null\quad\null& 
2n=a'r+\rho',\;0\les\rho'<r,
\\ && 
\text{and $x_1,\dots,x_{2n}\in\mbb P^1_{left}$ are distinct points.} 
\end{array}
$$

\nit(iii) In the case $r=2$, of mathematical instantons, Tikhomirov~\cite{tikh-odd,tikh-even} proved that their moduli spaces are irreducible, so the Egen-condition~\ref{def:FonP3}(iv) can be dropped. 
\end{m-remark}

%%%%%%%%%%%%%%%%%%%%%%%%
%%%%%%%%%%%%%%%%%%%%%%%%

\end{document}